\documentclass[11pt]{article}

\usepackage[a4paper]{geometry}
\usepackage[T1]{fontenc}
\usepackage{ae}
\usepackage{eucal}
\usepackage[utf8]{inputenc}
\usepackage{amsmath,amsfonts,amssymb}
\usepackage{url}
\usepackage{graphicx}

\newcommand{\legendre}[2]{\genfrac {(}{)}{1pt}{}{#1}{#2}}
\newcommand{\w}{\mathfrak w}

\newcommand{\OK}{\mathcal{O}_{K}}

\newcommand{\GFq}[1]{\mathbb{F}_{#1}}
\newcommand{\Z}{\mathbb Z}

\newcommand{\Q}{\mathbb Q}
\newcommand{\C}{\mathbb C}
\newcommand{\K}{\mathbb{K}}
\newcommand{\Kbf}{\mathbf{K}}

\newcommand{\wf}{\mathfrak{f}}
\newcommand{\wfone}{\mathfrak{f}_1}
\newcommand{\wftwo}{\mathfrak{f}_2}

\newtheorem{theorem}{Theorem}
\newtheorem{proposition}[theorem]{Proposition}

\newtheorem{lemma}[theorem]{Lemma}
\newtheorem{corollary}[theorem]{Corollary}

\newtheorem{example}[theorem]{Example}

\newenvironment {proofof}[1]
   {{\parindent0cm \bf Proof#1:}}
   {\hspace* {\fill} $\Box$}
\newenvironment {proof}
   {\begin {proofof}{}}
   {\end {proofof}}

\begin{document}

\title{Edwards curves and CM curves}

\author{
Fran\c{c}ois Morain \\
INRIA Saclay--\^Ile-de-France \\
\& Laboratoire d'Informatique (CNRS/UMR 7161) \\
\'Ecole polytechnique \\
91128 Palaiseau \\
France \\
morain@lix.polytechnique.fr}
\date{April 14, 2009}

\maketitle

\begin{abstract}
Edwards curves are a particular form of elliptic curves that admit a
fast, unified and complete addition law. Relations between Edwards
curves and Montgomery curves have already been described. Our work
takes the view of parameterizing elliptic curves given by their
$j$-invariant, a problematic that arises from using curves with
complex multiplication, for instance. We add to
the catalogue the links with Kubert parameterizations of $X_0(2)$ and
$X_0(4)$. We classify CM curves that admit an Edwards or Montgomery
form over a finite field, and justify the use of isogenous curves when
needed.
\end{abstract}


\section{Introduction}

An Edwards curve \cite{Edwards07} is a particular form of an elliptic
curve that leads to fast, unified and complete addition formulas
\cite{BeLa07,BeLa09}. Using
inverted Edwards coordinates yields the fastest formulas for
adding points on such a curve \cite{BeLa07b}. One may consult
\cite{BeLa07d} for comparisons of multiplication on a curve, as well
as for many references.

Following \cite{BeBiJoLaPe08}, the equation of a {\em twisted Edwards}
curve is 
\begin{equation}\label{eq:Edwards}
a x^2 + y^2 = 1 + d x^2 y^2
\end{equation}
where $a$, $d$ non zero elements of a field $\K$. An Edwards curve
corresponds to $a=1$.

The completeness of the addition
equations occurs only when $a$ is a square and $d$ is not a square in
the base field, in which case there are no rational singular points on
the curve. We will say that an elliptic curve $E$ can be
put in {\em twisted Edwards form} if there is a (rational) change of variables
leading to an equation of the type (\ref{eq:Edwards}). We will use the
term {\em complete Edwards form} when $a$ is a square (including the
case $a=1$ of course) and $d$ is not a square.

Recognizing an Edwards curve is relatively easy 
\cite[Theorem 3.3]{BeBiJoLaPe08}.
\begin{theorem}
If $\K$ is of characteristic different from $2$, the curve $E$ is
birationally equivalent to an Edwards curve if and only if $E(\K)$ has
a point of order 4.
\end{theorem}
By \cite[Theorem 2.1]{BeLa07}, if $\K$ is finite, the curve is
complete if $E(\K)$ has a point of order 4 and a unique point of order
2.

Another form of elliptic curve is the Montgomery form
\begin{equation}\label{eq:Mty}
E: B y^2 = x^3 + A x^2 + x
\end{equation}
with $A \neq \pm 2$, $B \neq 0$ (see \cite{Montgomery87}). Theorem 3.2 of
\cite{BeBiJoLaPe08} is
\begin{theorem}
If $\K$ has characteristic different from 2, then every curve in
Montgomery form is birationally equivalent to a twisted Edwards
curve, the converse being true.
\end{theorem}

In many applications (factoring, elliptic curve cryptography), one has
the choice of the curve and parameterization. In other contexts, such
as primality proving or the CM method, a curve is given via its
$j$-invariant and it may seem interesting to find an appropriate
parameterization starting from $j$.

Since elliptic curves having a rational $2$-torsion or $4$-torsion
point can be parameterized using modular curves, we will indicate how to
relate this to Edwards or Montgomery forms via Kubert's
equations. This will be the task of Section \ref{sct:2}, together with
more general properties of the $2$-torsion and $4$-torsion.
Section \ref{sct:CM} is interested in CM properties of elliptic
curves, and use classical results to investigate when the discriminant
of a curve having CM by a quadratic order is a square in the
corresponding ring class field. Results of Section
\ref{sct:2} and \ref{sct:CM} are then applied in Section \ref{sct:3}
on curves over finite fields. Moreover, we show how to replace a curve that
does not admit a complete Edwards form by an isogenous curve having such a
form, thus extending \cite[Section 5]{BeBiJoLaPe08}.

\section{Properties of the $2$-torsion}
\label{sct:2}

\subsection{Generalities}

We collect here some properties of the $2$-torsion group of an
elliptic curve whose equation will be taken as $Y^2 = F(X) = X^3 + a_2
X^2 + a_4 X + a_6$, with coefficients in some field $\K$ of
characteristic $\neq 2$ (see \cite{Silverman86}). We will note
$e_1$, $e_2$, $e_3$ for the roots of
$F(X)$ (over some subfield of an algebraic closure of the base field
$\K$). Remember that the discriminant of the curve is
$$\Delta(E) = 
-64\,a_6\,{a_2}^{3}+16\,{a_4}^{2}{a_2}^{2}
+288\,a_4\,a_6\,a_2-64\,{a_4}^{3}-432\,{a_6}^{2}.$$
Two curves $E: Y^2 = X^3 + a_2 X^2 + a_4 X + a_6$ and $E': Y^2 = X^3 +
a_2' X^2 + a_4' X + a_6'$ are isomorphic if and only if there exists
$(u, r)$ in $\K^2$ with $u\neq 0$ such that the following system has a
solution:
\begin{equation}\label{eq2}
\left\{\begin{array}{lcl}
u^2 a_2' &=& a_2 + 3 r, \\
u^4 a_4' &=& a_4 + 2 r a_2 + 3 r^2, \\
u^6 a_6' &=& a_6+r a_4+r^2 a_2+r^3,\\
\end{array}\right.
\end{equation}
the isomorphism sending $E$ to $E'$ via
$$(x, y) \mapsto (u^2 x+r, u^3 y).$$

The $n$-th division polynomial will be denoted $f_n(X)$ and we will be
interested mostly in the properties of
$$f_4(X) = 
{X}^{6}+2\,a_2\,{X}^{5}+5\,a_4\,{X}^{4}+20\,a_6\,{X}^{3}+ \left(
20\,a_2\,a_6-5\,{a_4}^{2} \right) {X}^{2}$$
$$\hspace*{3cm}+ \left(
-4\,a_4\,a_6+8\,{a_2}^{2}a_6-2\,a_2\,{a_4}^{2} \right)
X+4\,a_4\,a_6\,a_2-{a_4}^{3}-8\,{a_6}^{2}.$$
Note that the discriminant of $F(X)$ is $\Delta(E)/16$ and that of
$f_4(X)$ is $-\Delta(E)^5/4$.

\subsection{Curves with at least one torsion point}

We now turn our attention towards the $2$-torsion and $4$-torsion of
$E$. If $E$ has a rational $2$-torsion point, then $F(X)$ will have
one or three roots over $\K$. We will call these cases type I or type
III curves.

Let us begin with some general results on the curve $E$ of
equation $Y^2 = F(X) = (X-x_0) (X^2 + C X + D)$. Elementary
computations show the following list of results:
$$\Delta(E) = 16 \left( x_0^{2} + C x_0 +D \right) ^{2}
\left( {C}^{2} -4\,D \right),$$
which cannot be zero since $E$ is an elliptic curve. Furthermore
\begin{equation}
\label{eq1}
f_4(X) = \mathcal{P}_2(X) \mathcal{P}_4(X)
\end{equation}
where $\mathcal{P}_i$ has degree $i$ and
$$\mathcal{P}_2(X) = {X}^{2}-2\,x_0\, X-Cx_0-D,$$
$$\mathcal{P}_4(X) = {X}^{4}+2\,C\,{X}^{3}+6\,D\,{X}^{2} + \left(-8\,x_0
\, D +2\,C D +2\,x_0\,{C}^{2} \right)
X+(4 D -{C}^{2}) x_0^{2}+{D}^{2}.$$
The discriminants are
$$\mathrm{Disc}(\mathcal{P}_2) = 4 (x_0^2 + C x_0 + D) = 4 \mathcal{D}_2,$$
$$\mathrm{Disc}(\mathcal{P}_4) = -2^8 \, \left( {C}^{2}-4\,D \right) ^{3} 
\left(x_0 ^{2} + Cx_0 + D\right) ^{3} 
= -2^8 \left( {C}^{2}-4\,D \right) ^{3} \mathcal{D}_2^3.$$

\begin{lemma}\label{lem1}
The polynomials $\mathcal{P}_2$ and $\mathcal{P}_4$ do not have a
common root.
\end{lemma}
\begin{proof}
We compute
$$\mathrm{Resultant}_X(\mathcal{P}_2(X), \mathcal{P}_4(X)) 
= -16 (x_0^2 + C x_0 + D)^3 (C^2-4 D),$$
and both terms are non-zero, otherwise $F$ would have multiple roots.
\end{proof}

\subsubsection{Curves of type I}

Such a curve has equation $Y^2 = F(X) = (X-x_0) (X^2 + C X + D)$ with
the quadratic polynomial irreducible.
Let us study division by $2$ on $E$. Let $P = (x, y)$ be a rational
$4$-torsion point. Then $[2] P$ is the
$2$-torsion point $(x_0, 0)$. Writing the formulas for multiplication by
$2$, we get
$$\frac{F'(x)^2}{4 F(x)} - (C-x_0) - 2x = x_0$$
or $x$ is a root of
$$\mathcal{P}_{x_0}(X) = \left( {X}^{2}-2\,x_0\,X-D-x_0\,C\right) ^{2}
= \mathcal{P}_2(X)^2.$$

\begin{corollary}\label{coro1}
If $E$ is of type I, the polynomial $\mathcal{P}_4$ cannot have a
rational root.
\end{corollary}

\begin{proof}
Assume on the contrary that $\mathcal{P}_4$ has a rational root
$z$. Then $z$ would be sent by multiplication by 2 on the unique
rational $2$-torsion abscissa $x_0$. This would
imply $\mathcal{P}_{x_0}(z) = 0$, which by Lemma \ref{lem1} is impossible.
\end{proof}

\begin{proposition}\label{prop2}
The curve $E: Y^2 = (X-x_0) (X^2+ C X + D)$ of type I
has two rational $4$-torsion points if and only if
$\mathcal{D}_2$ is a square.
\end{proposition}

\begin{proof}
The rational roots of the polynomial $f_4$ are that of
$\mathcal{P}_2$ by Corollary \ref{coro1}.
Writing $z^2 = x_0^2 + C x_0 + D$, the polynomial $\mathcal{P}_2$ has
roots $x_{\pm} = x_0 \pm z$. Letting $y_{\pm}$ denote the ordinates, we
find that $y_{\pm}^2 = z^2 (C + 2 (x_0 \pm z))$. Since
$$(C + 2 (x_0 + z)) (C + 2 (x_0 -z)) = C^2 - 4 D$$
we see that exactly one of the factors is a square, leading to two
rational $4$-torsion points.
\end{proof}

\subsubsection{Curves of type III}
\label{sct:typeIII}

Suppose now $E: Y^2 = F(X) = (X-e_1) (X-e_2) (X-e_3)$ is of type III
with all $e_i\in\K$. Then
$$\Delta(E) = 2^4 (e_1-e_2)^2 (e_1-e_3)^2 (e_2-e_3)^2,$$
and the polynomial $\mathcal{P}_4$ will factor into three quadratics
$$f_4(X) = 
\left( {X}^{2}-2\,e_1\,X+e_1 (e_2+e_3)-e_2\,e_3 \right)$$
$$\hspace*{2cm}\times\left({X}^{2}-2\,e_2\,X+e_2 (e_1+e_3)-e_1\,e_3\right)$$
$$\hspace*{2cm}\times\left({X}^{2}-2\,e_3\,X+e_3 (e_1+e_2)-e_1\,e_2\right)$$
of respective discriminants
$4 (e_1 - e_3) (e_1 - e_2)$, $4 (e_2 - e_1) (e_2 - e_3)$, $4 (e_3-e_1)
(e_3-e_2)$. If $\Delta(E)$ is a square in $\K$, then one or all
these discriminants are squares, so that the corresponding factors
split. Suppose $4 (e_1 - e_3) (e_1 - e_2) = \delta_1^2$. Then the two roots of
$$\mathcal{Q}_{e_1}(X) = {X}^{2}-2\,e_1\,X+e_1 (e_2+e_3)-e_2\,e_3$$
are $x_{\pm} = 2 e_1 \pm \delta_1$. The corresponding ordinates satisfy
$$y_{\pm}^2 = \delta_1^2/4 (2 x_{\pm} - e_2 - e_3).$$ 

\subsection{Properties of $2$-isogenies}

The following formulas come
from the use of V\'elu's formulas (already given in \cite{CoDeMo96}):
\begin{proposition}\label{Velu2}
Assume that $E: Y^2 = X^3 + a_2 X^2 + a_4 X + a_6$ has a rational
point of order $2$, noted $P = (x_0, 0)$.
Put $t = 3 x_0^2 + 2 a_2 x_0 + a_4$ and $w = x_0 t$.
Then an equation of $E/\langle P\rangle$ is $E_1: Y_1^2 = X_1^3 +
A_2 X_1^2+ A_4 X_1 + A_6$ where $A_2 = a_2$, $A_4 = a_4 - 5 t$, $A_6 =
a_6 - 4 a_2 t - 7 w$. Moreover, the isogeny $I_1: E
\rightarrow E_1$ sends $(X, Y)$ to
$$(X_1, Y_1) = \left(X + \frac{t}{X-x_0}, Y
\left(1-\frac{t}{(X-x_0)^2}\right)\right).$$
\end{proposition}

The typical use of the Proposition \ref{prop2} is given now. It will
prove essential in the isogeny volcano approach later on.
\begin{corollary}\label{coro2}
Let $E$ be a type III curve. If $E_1 = E/\langle (e_1, 0)\rangle$ is a
type I curve, then $E_1$ admits a complete Edwards form.
\end{corollary}

\begin{proof}
An equation for $E_1 = E/\langle (e_1, 0)\rangle$ is
$$E_1: Y^2 = (X-x_0) (X^2 + C X + D)$$
with
$$x_0 = e_2 + e_3 - e_1, \quad C^2 - 4 D = 16 (e_1 - e_2) (e_1-e_3), 
\quad x_0^2 + C x_0 + D = (e_2-e_3)^2.$$
If $e_2$ and $e_3$ are rational, then $\mathcal{D}_2$ is a square and
we are done.
\end{proof}

\subsection{Montgomery parameterizations}

A Montgomery form of an elliptic curve is some equation $B y^2 = x^3 +
A x^2 + x$ with $A\neq \pm 2$ and $B\neq 0$. This shows that only curves
having a rational $2$-torsion point can be of this form.
Curves having a $2$-torsion points are points on the modular curve
$X_0(2)$, an equation of which is
\begin{equation}\label{eq3}
j = \frac{(u+16)^3}{u}.
\end{equation}
Note also that
$$j - 1728 = \frac{(u+64) (u-8)^2}{u}$$
An equation for an elliptic curve $E$ of given invariant $j$ is
\begin{equation}\label{eq4}
Y^2 = X^3 + \frac{3 j}{1728-j} c^2 X + \frac{2 j}{1728-j} c^3
\end{equation}
(with $c$ present to accommodate twists).
We compute
$$\Delta(E) = 2^{12} 3^6 c^6 \frac{j^2}{(j-1728)^3} = 
2^{12} 3^6 c^6\,{\frac { \left( u+16 \right) ^{6}u}{ \left( u+64 \right) ^{3}
 \left( u-8 \right) ^{6}}}.$$
Plugging equation (\ref{eq3}) in (\ref{eq4}), we find that the cubic
has rational root $- c (u+16)/(u-8)$. The quadratic factor has
discriminant $\mathcal{D}_1 = 9 c^2 u (u+64) (u+16)^2$. Moreover
\begin{equation}\label{eq6}
\mathcal{D}_2 = 2^6\times 3^2 c^2 (u+16)^2 (u-8)^2 (u+64).
\end{equation}

Let us show how we can recover a Montgomery parameterization for these
curves. Setting $X = X' - c \frac{u+16}{u-8}$, this transforms
(\ref{eq4}) into
$$Y^2 = {X'}^3 -3\,c\,{\frac { \left( u+16 \right) {{X'}}^{2}}{u-8}}+
144\,c^2\,{\frac { \left( u+16 \right) ^{2}{X'}}{ \left( u+64 \right) 
 \left( u-8 \right) ^{2}}}.$$
If $u+64 = v^2$, then setting $X' = k X''$ with
$$k = 12\,c\,{\frac {u+16}{ \left( u-8 \right) v}}$$
leads to the Montgomery form
$$\frac{1}{k} (Y/k)^2 = {X''}^3 - \frac{v}{4} {X''}^2 + {X''}.$$
In Section \ref{sct:CM}, we will express $u$ as one of Weber's
functions and investigate when $u+64$ is a square in some ring class
field.

\subsection{Kubert parameterizations}
\label{ssct:Kubert}

In \cite{Kubert76} are given all parameterizations of curves over $\Q$
having prescribed torsion structure. Of particular relevance to
Edwards curves is the parameterization for curves $E$ containing the
torsion group $\Z/4\Z$. Letting $b\in\Q$ such that $b^4 (1+16 b)\neq 0$, we
get the parameterization
$$\mathcal{EK}_b: Y^2 = (X-4 b)  \left( {X}^{2}+ X - 4 b \right).$$
The curve $\mathcal{EK}_b$ has a unique point of order $2$, $(4 b, 0,
1)$. The division polynomial $f_4$ factors as
$$f_4(X) = X \left( X-8\,b \right)
\left({X}^{4}+2\,{X}^{3}-24\,b{X}^{2}+128\,{ b}^{2}X-256\,{b}^{3}
\right).$$
The two rational roots are: $0$, which leads to two
rational points $(0, \pm 4 b, 1)$ and $8b$, for which $Y^2 =
16 b^2(16b+1)$.
Note that the $j$-invariant of $\mathcal{EK}_b$ is:
$$j = {\frac { \left( 16\,{b}^{2}+16\,b+1 \right) ^{3}}{{b}^{4} \left(
16\,b +1 \right) }}.$$
Writing $w = 1/b$ leads to
$$j = {\frac { \left( {w}^{2}+16\,w+16 \right) ^{3}}{w \left( 16+w
\right) }}$$
and we see that setting $u = w^2+16 w$ makes $j$ of the form (\ref{eq3})
and that $u + 64 = (w+8)^2$, so that we get a Montgomery form in that
case too. With the notations of the preceding Section:
$$\mathcal{D}_1 = 9 c^2 w^2 (w+8)^2 (w^2+16 w+16)^2 (1+16/w),$$
$$\mathcal{D}_2 = 2^6\times 3^2 c^2 (w^2+16 w+16)^2 (w^2+16 w-8)^2 (w+8)^2.$$
This shows the following
\begin{proposition}
For any $w\neq 0$, the curve $E$ associated to $w$ admits a Montgomery
form. Moreover, if $1+16/w$ is not a square, $E$ admits a complete
Edwards form.
\end{proposition}

Making $w = 16 d / (a-d)$ or $a = d (1 + 16/w)$ yields directly 
$$J = \frac{16 (a^2 + 14 a d + d^2)^3}{a d (a-d)^4}$$
and this is precisely the invariant of a twisted Edwards curve.

The proof of the following is rather tedious and is preferably done
using a computer (and {\sc Maple} in our case).
\begin{proposition}
Suppose $E$ is of type I: $F(X) = (X-x_0)(X^2+C X+D)$ with the
quadratic polynomial irreducible and $x_0^2 + C x_0 + D = z^2$. Then
$E$ is isomorphic to $\mathcal{EK}_b$ with
$$r = \frac{(1-4b) \,{u}^{2} + x_0 -C}{3},$$
$$u^2 = C + 2 (x_0 \pm z),$$
whichever sign yields a square (cf. Proposition \ref{prop2}) and
$$b = \mp 1/4\,{\frac {z}{C+2\, (x_0 \pm \,z)}}
\text{ or }
b = -1/4\,{\frac {\pm 7\,z+8\,x_0+4\,C}{C+2\, (x_0\pm \,z)}}.$$
\end{proposition}

\begin{proposition}
The Kubert curve $\mathcal{EK}_b$ is birationally equivalent to the
Edwards curve of parameter $d = 16 b+1$. If $16b+1$ is not a square,
then the curve is complete.
\end{proposition}

\begin{proof}
If we plug these values in the transformation formulae of
\cite[Theorem 2.1]{BeLa07}, we get
$$a_2 = 2 b+1/4, a_4 = b^2, d = 16 b+1, r_1/(1-d) = 1/16$$
so that the curve $\mathcal{EK}_b$ is birationally equivalent to
$$E': (s/4)^2 = r^3+a_2r^2+a_4r,$$
which is shown in the same theorem to be birationally equivalent to
$x^2+y^2 = 1 + d x^2 y^2$.
\end{proof}

\section{CM curves}
\label{sct:CM}

In the CM method, we use elliptic curves that are constructed given
their invariant. We have the choice of the explicit form of
the equation to be used and it is natural to ask when a CM curve can
be written in Edwards or Montgomery form (see \cite{Bernstein06} for a
possible use in primality proving \cite{AtMo93b}).

\subsection{Theorems over $\C$}

Let $\mathcal{E}$ have complex multiplication by an order $\mathcal{O}$ of
discriminant $D = f^2 D_K$ in an
imaginary quadratic field $\Kbf = \Q(\sqrt{D_K})$ of discriminant
$D_K$. Such a curve can be built using its $j$-invariant that is a
root of the so-called {\em class polynomial} $H_D(X)$, that generates the
{\em ring class field} $\Kbf_{\mathcal{O}} =
\Kbf(j(\mathcal{O}))$. The roots of $H_D(X)$ are of the form
$j(\alpha)$ where $\alpha$ is the root of positive imaginary part of
$A X^2 + B X+ C$ with $(A, B, C)$ a reduced primitive quadratic form
of discriminant $B^2-4A C = D = \mathrm{Disc}(\alpha)$.

Given $j(\alpha)$ (in other words, any root of $H_D(X)$),
we can use equation (\ref{eq4}) for $\mathcal{E}(\alpha)$ and
look for cases where $\Delta(\mathcal{E}(\alpha))$
is a square in $\Kbf_{\mathcal{O}}$. These results will translate in
equivalent properties over finite fields.
It is natural for this to introduce the Weber function $\gamma_3$
satisfying $\gamma_3(\alpha)^2 = j(\alpha)-1728$. We rephrase 
\cite[Satz (5.2)]{Schertz76} (see also \cite{Schertz02}) as:
\begin{theorem}\label{thm0}
(a) When $D$ is odd, $\Delta(\mathcal{E}(\alpha))$ is a square in $\Q(\alpha,
j(\alpha))$.

(b) When $D$ is even, $\Delta(\mathcal{E}(\alpha))$ is a square in
$\Q(j(2 \alpha))$.
\end{theorem}

\begin{proof}
(a) When $D$ is odd $\sqrt{D}\gamma_3(\alpha)$ is a class invariant
and we write
$$j(\alpha)-1728 = \left(\frac{\sqrt{D}\gamma_3(\alpha)}{\sqrt{D}}\right)^2.$$

(b) We use $\Q(\gamma_3(\alpha)) = \Q(j(2 \alpha))$.
\end{proof}

To understand when CM curves have rational $2$-torsion points, we use
some other Weber functions, namely $\wf(\alpha)$, $\wfone(\alpha)$ and
$\wftwo(\alpha)$ that satisfy
$$j(\alpha) =
\frac{(-\wf(\alpha)^{24}+16)^3}{-\wf(\alpha)^{24}} =
\frac{(\wfone(\alpha)^{24}+16)^3}{\wfone(\alpha)^{24}} =
\frac{(\wftwo(\alpha)^{24}+16)^3}{\wftwo(\alpha)^{24}},$$
in other words, $-\wf(\alpha)^{24}$, $\wfone(\alpha)^{24}$ and
$\wftwo(\alpha)^{24}$ are the roots of (\ref{eq3}).
The numbers $\wf(\alpha)^{24}$, $\wfone(\alpha)^{24}$ or
$\wftwo(\alpha)^{24}$ are very often {\em class invariants}, that is
elements of $\Q(j(\alpha))$; sometimes they are in $\Q(\alpha,
j(\alpha))$. Moreover small powers of these functions are very often
elements of $\Q(j(r \alpha))$ for some $r = 2^{\pm n}$, as can be seen
from \cite{Schertz76} for instance.

To go further, we introduce the generalized Weber function
$$\w_N(z)^s = (\eta(z/N)/\eta(z))^s$$
where $\eta$ is Dedekind's function and
$N$ an integer and $s$ some integer related to $N$. These functions
are modular for $\Gamma^0(N)$ and give a model for $X^0(N)$ (equivalently
$X_0(N)$). In particular \cite{EnMo09}
\begin{theorem}
Let $\alpha = \frac{-B+\sqrt{D}}{2A}$ be a root associated to the
primitive reduced quadratic form $[A, B, C]$ of discriminant $D = B^2
- 4 A C$. If $B^2 \equiv D \bmod 16$ has a solution in $B$
(equivalently $D\bmod 16 \in \{0, 1, 4, 9\}$), then 
$\Q(j(\alpha)) \subset \Q(\w_4^8(\alpha)) \subset \Q(\alpha, j(\alpha))$.
\end{theorem}
From this
\begin{corollary}\label{coro3}
Let $D\bmod 16 \in \{0, 1, 4, 9\}$. Then $\mathcal{E}(j(\alpha))$
admits a Montgomery form.
\end{corollary}
\begin{proof}
Use the fact that the modular equation linking $j(z)$ and $\w_4(z)$ is
precisely
$$J = {\frac { \left( {w}^{2}+16\,w+16 \right) ^{3}}{w \left( 16+w
\right) }}$$
which sends us back to Section \ref{ssct:Kubert}.
\end{proof}

\section{Properties of the $2$-torsion over prime finite fields}
\label{sct:3}

\subsection{Splitting properties}

Part of what follows can also be found in \cite{Sutherland08}.
Let $\K$ be a prime finite field of characteristic $p>2$.
The following result is classical and taken from \cite{Swan62}. It
will help us study the splitting properties of $F(X)$ and
$f_4(X)$ over a finite field. Note that $F(X)$ and $f_4(X)$ have no
square factor (since $\Delta(E)\neq 0$ for $E$ to be an elliptic curve).
\begin{theorem}
Let $f(X)$ be a squarefree polynomial of degree $d$ and $n$ its number of
irreducible factors modulo $p$. Then
$$\legendre{\mathrm{Disc}(f)}{p}  = (-1)^{d-n}.$$
\end{theorem}
This gives us immediately.
\begin{proposition}\label{prop3}
Let $p$ be an odd prime.
The curve $E$ has exactly one $2$-torsion point over $\GFq{p}$ if and
only if $\legendre{\Delta(E)}{p} = -1$. In that case, and writing $n_4$
for the number of irreducible factors of $f_4$, one has
$$(-1)^{n_4} = -\legendre{-\mathcal{D}_2}{p}.$$
\end{proposition}
\begin{proof}
We have
$$\legendre{\mathrm{Disc}(\mathcal{P}_4)}{p} =
\legendre{-1}{p}\legendre{{C}^{2}-4\,D}{p}\legendre{\mathcal{D}_2}{p}.$$
If $E$ has a unique $2$-torsion point, then
$$\legendre{{C}^{2}-4\,D}{p}=-1$$
which yields the result.
\end{proof}

When a polynomial $P(X)$ has factors of degrees $d_1$, $\ldots$,
$d_k$ over $\K$, we will denote this splitting as $(d_1) \cdots (d_k)$.
The following Proposition describes the splittings of $\mathcal{P}_2$
and $\mathcal{P}_4$.
\begin{proposition}
Let $E: Y^2 = F(X) = (X-x_0) (X^2+ C X +D)$ be of type I.
The splittings of $\mathcal{P}_2$ and $\mathcal{P}_4$ can be found in
the following table:
$$\begin{array}{|c||c|c||c|c|}\hline
                 & \multicolumn{2}{c||}{\legendre{\mathcal{D}_2}{p} = +1}
		 & \multicolumn{2}{c|}{\legendre{\mathcal{D}_2}{p} = -1}
 \\ \hline
                 & \mathcal{P}_2 & \mathcal{P}_4 & \mathcal{P}_2 &
		 \mathcal{P}_4 \\ \hline
p\equiv 1\bmod 4 & (1) (1) & (4) & (2) & (2) (2) \\
p\equiv 3\bmod 4 & (1) (1) & (2) (2) & (2) & (4) \\
\hline
\end{array}$$
\end{proposition}

\begin{proof}
Assume first that $\legendre{\mathcal{D}_2}{p} = +1$.
If $p\equiv 1\bmod 4$, $\mathcal{P}_4$
should have an odd number of irreducible factors, leading to $(4)$ or
$(1)(1)(2)$ but Corollary \ref{coro1} rules out $(1)(1)(2)$. If $p\equiv
3\bmod 4$, then $\mathcal{P}_4$ should have an even number of factors
or be of type $(1)(3)$ and $(2)(2)$ and only the latter one is
possible.

The proof for the case $\legendre{\mathcal{D}_2}{p} = -1$ is
symmetrical and we omit it.
\end{proof}

\subsection{Reduction of CM curves over a finite field}

If $p$ splits in the
ring class field $\Kbf_{\mathcal{O}} = \Kbf(j(\mathcal{O}))$, i.e., $p =
(U^2 - D V^2)/4$, then we can reduce $\mathcal{E}$ modulo a prime
factor of $p$ in $\Kbf_{\mathcal{O}}$ to get a curve $E/\GFq{p}$ of
cardinality $p+1-U$. This is the heart of the CM-method, which is a
building block in ECPP for instance \cite{AtMo93b,Morain07}.

Conversely, a (non supersingular) elliptic curve $E/\GFq{p}$ has complex
multiplication by an order $\mathcal{O}$ in an imaginary quadratic
field $\Kbf$. In details, if $E$ has cardinality $p+1-U$, write
$\mathrm{Disc}(\pi) = U^2-4 p = V^2 D_K$
for the discriminant of the Frobenius $\pi$ of the curve.
We have $\Z[\pi] \subset \mathcal{O} \subset \OK$, where $\OK$
is the ring of integers of $\Kbf$.
Noting $f$ for the conductor of the order $\mathcal{O}$, we have
$\mathrm{Disc}(\mathcal{O}) = f^2 D_K$, with $f \mid V$.

\begin{figure}[Hbt]
\iftrue
\centerline{\includegraphics{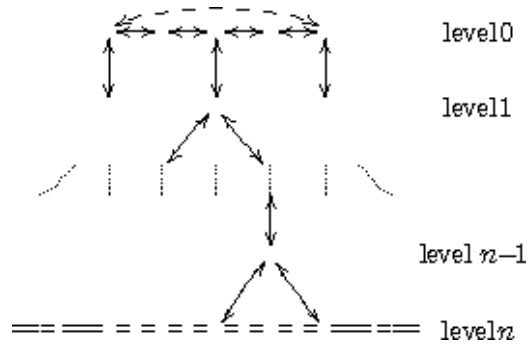}}
\else
\begin{center}
$\xymatrix@-0.8pc{
& 
 & 
 & \ar@{<-->}@/^1pc/ [rrrr] \ar@{<->}[d] \ar@{<->}[r] 
 & \ar@{<->}[r] 
 & \ar@{<->}[d] \ar@{<->}[r]
 & \ar@{<->}[r] 
 & \ar@{<->}[d]
 &
 &
 &
 & \text{level } 0\\
& 
 & 
 & 
 &
 & \ar@{<->}[dr] \ar@{<->}[ld]
 & 
 & 
 & 
 & 
 &
 & \text{level } 1\\
& & \ar@{..}[dl] 
  & \ar@{..}[d] 
  & \ar@{..}[d] 
  & \ar@{..}[d] 
  & \ar@{..}[d]
  & \ar@{..}[d]
  & \ar@{..}[dr] 
  &
  & \\ 
& 
 & 
 & 
 &
 & 
 & \ar@{<->}[d]
 & 
 & 
 & 
 & 
 & \\
& 
 & 
 & 
 &
 & 
 & \ar@{<->}[dr] \ar@{<->}[ld]
 & 
 & 
 & 
 & 
 & \text{level } n-1\\
 \ar@{==}[r] & \ar@{==}[r] & \ar@{==}[r] &
 \ar@{==}[r] & \ar@{==}[r] & \ar@{==}[r] & \ar@{==}[r] & \ar@{==}[r] &
 \ar@{==}[r] & \ar@{==}[r] & & \text{level } n}$ 
\end{center} 
\fi
\caption{Isogeny volcano for the prime $2$.}  \label{volcano} 
\end{figure}

The volcano structure \cite{Kohel96} (see also \cite{FoMo02})
describes the relationships between inclusions
of orders in $\OK$ and the structure of elliptic curves having CM by
these orders. Rational $2$-torsion points are in one-to-one
correspondence with isogenies of degree $2$ (see Proposition
\ref{Velu2} for an illustration of this). The volcano for the prime 2
has the shape of Figure \ref{volcano} (in case $(2)$ splits in
$\OK$). The crater is formed of
horizontal isogenies (if any) and each curve on the crater has one
isogeny down. Any curve strictly between the crater and the bottom has
one isogeny up and two down.

General properties of the
volcano can be used to meet some of our needs. The following results
justifies the idea already presented in \cite[Section 5]{BeBiJoLaPe08}.
\begin{theorem}\label{propfdale}
Assume $E/\GFq{p}$ is of type III. There exists a curve $E'/\GFq{p}$
isogenous to $E$ that is of type I. Moreover, $E'$ admits a complete
Edwards form.
\end{theorem}

\begin{proof}
The proof of Proposition \ref{propfdale} comes directly
from \cite[Proposition 23]{Kohel96}.
It is enough to find a curve at the bottom of the volcano for prime $2$.
Letting $\mathrm{Disc}(\pi) = U^2-4p = D_K V^2$ and $2^n \mid\mid V$,
this curve will be at level $n$.

The last part comes from Corollary \ref{coro2}.
\end{proof}

From a practical point of view, procedure {\sc FindDescendingPath} of
\cite{FoMo02} can be used to do just this. Now we apply the
results of Section \ref{sct:2}. Starting from our curve $E$, we build
a path
$$E \rightarrow E_1 \rightarrow \cdots \rightarrow E_{n-1} \rightarrow
E'$$
with $E_{n-1}$ of type III and $E'$ of type I. Proposition \ref{Velu2}
can be used to transport points if needed.

\medskip
\noindent
{\bf Remarks.} 

1. We can easily modify procedure {\sc FindDescendingPath}
so that we can keep track of the $2$-torsion points we encounter, so
that we obtain a root of $E/\langle (e_1, 0)\rangle$ from the $e_i$'s,
see Corollary \ref{coro2}.

2. If we know that $E$ has CM by $\OK$, then we can simplify {\sc
FindDescendingPath} by discarding horizontal invariants using the
class polynomial, so that we are left with just one descending path.

\medskip
\noindent
{\bf Numerical examples.} Consider $E: Y^2 = X^3 + X + 2$ over
$\GFq{1009}$. We find
$$\begin{array}{|c|c|c|}\hline
i & E_i & 2\text{-torsion} \\ \hline
0 & {[1, 2]} & \{463, 547, 1008\} \\
1 & {[990, 30]} = E_0/\langle (1008, 0)\rangle & \{2, 3, 1004\} \\
2 & {[950, 871]} = E_1/\langle (3, 0)\rangle & \{265, 750, 1003\} \\
3 & {[1003, 17]} = E_2/\langle (750, 0)\rangle & \{518\} \\
\hline
\end{array}$$
and $E_3$ admits $(247, \pm 19)$ has rational $4$-torsion points, and
therefore is birationally equivalent to a complete Edwards curve.

\medskip
We end this section with the following count, that is easily deduced
from the volcano structure.
\begin{proposition}\label{prop4}
The set of invariants of complete Edwards curves is formed of the
$j$-invariants on the floor of the $2$-volcano. If $4p = t^2 - v^2 2^n
D_K$, where $D_K$ is fundamental and $v$ odd, there is a total of 
$2^{n-1} (2 - \legendre{D_K}{2}) h(D_K)$ such invariants, where
$\legendre{D_K}{2}$ denotes the Kronecker symbol.
\end{proposition}

\subsection{Classifying Montgomery and Edwards curves over finite fields}

We assume throughout that $E/\GFq{p}$ is the reduction of a curve
$\mathcal{E}(\alpha)$ having CM by an order of
discriminant $D$ so that in particular $p = (U^2 - D V^2)/4$.
The aim of this section is to prove the following results.
\begin{theorem}\label{thm2}
If $D$ is fundamental, then $E$ does not admit a complete Edwards form.
\end{theorem}
\begin{theorem}\label{thm4}
Suppose $E/\GFq{p}$ has CM by $\mathcal{O}$ of discriminant $D$. The
following Table summarizes the properties of the reduction of
$\mathcal{E}(\alpha)$:
\begin{center}
\begin{tabular}{|c|c|c|c|c|}\hline
$D$       & $V$ & 2-torsion & Montgomery form & Edwards form \\ \hline
\multicolumn{5}{|c|}{$D$ odd} \\ \hline
$1 \bmod 8$ & --  & type III  & yes & twisted, not complete \\ \hline
$5 \bmod 8$ & even& type III  & yes & twisted, not complete \\
            & odd & none & -- & -- \\ \hline
\multicolumn{5}{|c|}{$D$ even} \\ \hline
$0, 4 \bmod 16$ & even & type III & yes & twisted, not complete \\
              & odd & type I & yes & complete \\ \hline
$8, 12 \bmod 16$ & even & type III & yes/no & twisted at best \\
                 & odd & type I & no & no \\ \hline
\end{tabular}
\end{center}
\end{theorem}

The following result is taken from \cite{Morain07b} (precised by
Theorem \ref{thm0}) and starts our proof of Theorem \ref{thm2}.
\begin{proposition}\label{prop1}
The quantity
$\Delta(E)$ is a square modulo $p = (U^2 - D V^2)/4$ in the following cases:

(a) $D$ odd;

(b) $D$ even and $2\mid V$.
\end{proposition}
Together with Proposition \ref{prop3}, this proves part of our
theorem.

\begin{proof}
By Section \ref{sct:typeIII}, any curve of type III admits rational
roots for $f_4(X)$, but not always rational ordinates.

Suppose $D$ satisfies one of the conditions of Proposition \ref{prop1}.
In that case, we see that $E$ has zero or three $2$-torsion
points. If we prove that $E$ admits at least one, we get three of them.

If $D\equiv 1\bmod 8$, the ideal $(2)$ splits in $\OK$ and there are three
(distinct) rational isogenies starting from $E$, corresponding to
three $2$-torsion points (\cite[Proposition 23]{Kohel96} again). Then
apply Corollary \ref{coro3}.

It is clear that when $D\equiv 5\bmod 8$ and $V$ odd,
then $E$ has no rational $2$-torsion points at all. 
With $p = (U^2-D V^2)/4$, we must have $U$ and $V$ of the same
parity. Having a $2$-torsion point is equivalent with $U$ even and
therefore $V$ even.

Suppose now $D = -4 m$ is even. Then $U$ is even, so that $E$ admits
at least one $2$-torsion point. If $V$ is even, $\Delta$ is a square,
forcing the splitting $(1)(1)(1)$. The cases $m = 0, 3 \bmod 4$ come
from Corollary \ref{coro3}.

Dirichlet's theorem (see \cite[Ch. 4]{Buell89}), gives us necessary
arithmetical conditions on $p$ to split as $(U^2 - D V^2)/4$.
For an integer $p$, let $\chi_4(p) = \legendre{-1}{p}$ and
$\chi_8(p) = \legendre{2}{p}$. The {\em generic characters} of $D$
are defined as follows:
\begin{itemize}
 \item $\legendre{p}{q}$ for all odd primes $q$ dividing $D$;
 \item if $D$ is even:
  \begin{itemize}
   \item $\chi_4(p)$ if $D/4 \equiv 3, 4, 7\pmod 8$;
   \item $\chi_8(p)$ if $D/4 \equiv 2\pmod 8$;
   \item $\chi_4(p) \cdot \chi_8(p)$ if $D/4 \equiv 6\pmod 8$;
   \item $\chi_4(p)$ and $\chi_8(p)$ if $D/4 \equiv 0\pmod 8$.
  \end{itemize}
\end{itemize}
\begin{theorem}\label{thm-Dirichlet}
An integer $p$ such that $\gcd(p, 2cD)=1$ is representable by some
class of forms in the principal genus of discriminant $D$ if and only
if all generic characters $\chi(p)$ have value $+1$.
\end{theorem}

The following finishes the proof of Theorem \ref{thm2}.
Suppose $D = -4 m$, $V$ odd and $m \bmod 4 \in \{1, 2\}$. We now show
that $E$ does not admit a $4$-torsion point.
If $m$ is odd, the equation $p = (U/2)^2 - (D/4) V^2 = (U/2)^2 + m V^2$
shows that $U/2$ should be even and therefore $p+1-U\equiv 2\bmod
4$, since $p\equiv 1\bmod 4$ as $\chi_4(p) = +1$.

Let $m$ be even. This implies $U \equiv 2\bmod 4$.
Suppose $m \equiv 2\bmod 8$. Write $p \equiv (U/2)^2 + 2 V^2 \bmod 8$.
Since $U/2$ must be
odd, we get $p \equiv 1 + 2 V^2\bmod 8$. On the other hand, $\chi_4(p)
\cdot \chi_8(p) = +1$, leading to $p\equiv 1\bmod 8$ and $V$ even or
$p\equiv 3\bmod 8$ and $V$ odd. In the latter case, we have $p+1-U
\equiv 4 - 2 \equiv 2\bmod 4$ and no rational $4$-torsion exists.

When $m \equiv 6\bmod 8$, we get $1 + 6 V^2 \equiv \pm 1\bmod 8$,
and $\chi_8(p)=+1$ implies $p\equiv 1\bmod 8$ and $V$ even; or
$p\equiv 7\bmod 8$, $V$ odd and $p+1-U \equiv 2\bmod 4$.

\medskip
We are left with the case $m \equiv 0\bmod 4$ and $m \equiv 3\bmod 4$,
in which case $D$ is not fundamental (or the order is not the
principal one). Assume $V$ odd.
Again using isogenies, we see that only one isogeny up goes from
$E$, so that $E$ admits only one rational
$2$-torsion point. Now apply Corollary \ref{coro3}.
\end{proof}

\subsection{Using isogenous curves}

As already stated, we can use some isogeny to get an Edwards curve
whenever possible in the same isogeny class as a given curve $E$. In
the case $D_K\equiv 1\bmod 8$ and $n=1$, we can also compute directly
a curve of type I having CM by $2 \OK$, since $h(4D_K) = h(D_K)$.

\section{Conclusions}

We have shed some light on the links between different
parameterizations and the Montgomery and Edwards form of an elliptic
curve. Curves with CM by a principal order cannot be of complete
Edwards form, though they may admit a Montgomery parameterization.
In practice, say in the course of the CM method \cite{AtMo93b}, this
could appear as a problem, but in many applications, replacing a curve by
an isogenous one is no harm.

\paragraph {Acknowledgments.} The author wants to thank the University
of Waterloo for its hospitality during his sabbatical leave. Thanks
also to G.~Bisson for directing my attention to \cite{BeBiJoLaPe08}
and useful discussions, and to A.~Sutherland for exchanging ideas on
the subject and for Proposition \ref{prop4}.

\iffalse
\bibliographystyle{plain}
\bibliography{morain}
\else
\def\noopsort#1{}\ifx\bibfrench\undefined\def\biling#1#2{#1}\else\def\biling#1%
#2{#2}\fi\def\Inpreparation{\biling{In preparation}{en
  pr{\'e}paration}}\def\Preprint{\biling{Preprint}{pr{\'e}version}}\def\Draft{%
\biling{Draft}{Manuscrit}}\def\Toappear{\biling{To appear}{\`A para\^\i
  tre}}\def\Inpress{\biling{In press}{Sous presse}}\def\Seealso{\biling{See
  also}{Voir {\'e}galement}}\def\Editor{\biling{Ed.}{R{\'e}d.}}

\fi

\end{document}